# Seven Conjectures in Geometry and Number Theory

Florentin Smarandache, Ph D
Professor of Mathematics
Chair of Department of Math & Sciences
University of New Mexico
200 College Road
Gallup, NM 87301, USA
E-mail: smarand@unm.edu

**Abstract:**
In this short paper we propose four conjectures in synthetic geometry that generalize Erdos-Mordell Theorem, and three conjectures in number theory that generalize Fermat Numbers.

**2000 MSC:** 11A41, 51F20

1. **Four Geometrical Conjectures:**

   a) Let $M$ be an interior point in a $A_1 A_2 ... A_n$ convex polygon and $P_i$ the projection of $M$ on $A_i A_{i+1}$, $i = 1, 2, 3, ..., n$.
   Then
   $$\sum_{i=1}^{n} \overline{MA_i} \geq c \sum_{i=1}^{n} \overline{MP_i}$$
   where $c$ is a constant to be found.

   For n=3, it was conjectured by Erdös in 1935, and solved by Mordell in 1937 and Kazarinoff in 1945. In this case $c = 2$ and the result is called the Erdös-Mordell Theorem.

   b) More generally: If the projections $P_i$ are considered under a given oriented angle $\alpha \neq 90$ degrees, what happens with the above inequality?

   c) In a 3-space, we make the same generalization for a convex polyhedron with n vertexes and m faces:
   $$\sum_{i=1}^{n} \overline{MA_i} \geq c_1 \sum_{j=1}^{m} \overline{MP_j}$$
   where $P_j$, $1 \leq j \leq m$, are projections of $M$ on all faces of the polyhedron, and $c_1$ is a constant to be determined.

[Kazarinoff conjectured that for the tetrahedron

$$\sum_{i=1}^{4} \overline{MA_i} \geq 2\sqrt{2} \sum_{i=1}^{4} \overline{MP_i}$$

and this is the best possible].

d) Furthermore, does the below inequality hold?

$$\sum_{i=1}^{n} \overline{MA_i} \geq c_2 \sum_{k=1}^{r} \overline{MT_k}$$

where $T_k$, $1 \leq k \leq r$, are projections of $M$ on all sides of the polyhedron, and $c_2$ is a constant to be determined.

## 2. Three Number Theory Conjectures (Generalization of Fermat Numbers):

Let's consider $a$, $b$ integers $\geq 2$ and $c$ an integer such that $(a, c) = 1$.

One constructs the function $P(k) = a^{b^k} + c$, where $k \in \{0, 1, 2, ...\}$.

Then:

a) For any given triplet $(a, b, c)$ there is at least a $k_0$ such that $P(k_0)$ is prime.

b) Does there exist a non-trivial triplet $(a, b, c)$ such that $P(k)$ is prime for all $k \geq 0$?

c) Is it possible to find a triplet $(a, b, c)$ such that $P(k)$ is prime for infinitely many $k$'s?

### REFERENCES

[1] Alain Bouvier and Michel George, sous la direction de François Le Lionnais, *Dictionnaire des Mathématiques Elémentaires*, Presses Universitaires de France, Paris, 1979.
[2] P. Erdös, Letter to T. Yau, August 1995.
[3] Florentin Smarandache, *Collected Papers, Vol. II*, University of Chişinău Press, Chişinău, 1997.

*[Published in author's book Collected Papers, Vol. II, 1997.]*